\documentclass[12pt]{article}
\pdfoutput=1
\usepackage{latexsym}
\usepackage{enumerate}
\usepackage{epsf,epsfig,amsfonts,a4wide}
\usepackage{mathtools}
\usepackage{amsfonts}
\usepackage{url}
\usepackage{amsmath,amssymb,amsthm}
\usepackage{epsf,epsfig,amsfonts,graphicx}
\usepackage{a4wide}

\usepackage[centerlast]{caption}
\usepackage{graphics,amsmath,amssymb,amscd}
\usepackage{amsbsy}
\usepackage{amsthm}
\usepackage{gensymb}
\usepackage{float}
\usepackage{graphicx}
\usepackage[colorlinks=true,linkcolor=black,citecolor=black]{hyperref}
\usepackage{multicol}
\usepackage{color}

\usepackage{epstopdf}

\newcommand{\be}{\begin{equation}}
\newcommand{\ee}{\end{equation}}

\usepackage{amsthm}
\usepackage{mathdots}

\usepackage{amscd}

\newcommand{\eps}{\varepsilon}
\newcommand{\ph}{\varphi}

\title{Bifurcations of phase portraits of pendulum with vibrating suspension point}
\author{A.I. Neishtadt$^{1,2,}$\footnote{Corresponding author.  E-mail addresses:  {\it A.Neishtadt@lboro.ac.uk} (A.I.Neishtadt),  {\it shengkaicheng@yeah.net} (K.Sheng)}
, K. Sheng$^1$\\ 
$^1$ Loughborough University, Loughborough, LE11 3TU, UK\\
$^2$ Space Research Institute, Moscow, 117997, Russia}

\date{}
\begin{document}
\maketitle

\begin{abstract}
We consider a simple pendulum whose suspension point undergoes fast vibrations in the plane of motion of the pendulum. The  averaged over the fast vibrations system is a   Hamiltonian system with one degree of freedom depending on two parameters. We give complete description of bifurcations of phase portraits of this  averaged system.

\end{abstract}

\vskip 20pt

\section{Introduction}

 A simple pendulum with vibrating suspension point is a classical problem of perturbation theory.  The phenomenon of stabilisation  of  the upper vertical position of the pendulum by fast vertical vibrations of the suspension point  was discovered by A. Stephenson \cite{steph_1, steph_2}.  In these papers the linearisation and reduction to the Mathieu equation is used.  The case of inclined vibrations of  the suspension point is considered as well.     Nonlinear theory was developed by N.N.Bogolyubov \cite{bogol_1}, who used the averaging method, and by P.L.Kapitsa, who has developed a method of separation of slow and fast motions for this  \cite{kapitsa_1, kapitsa_2} (see also \cite{LL}). Different aspects of this problem were discussed in many publications (see, e.g., \cite{levi} for a discussion of geometric aspects, \cite{bam} for the case of arbitrary frequencies of vibrations, \cite{aio} for the case of random vibrations). Generalisations to  double- and multiple-link pendulums are contained in \cite{steph_3, aches, khol}. It is noted in \cite{burd} that the problem is  simplified by using averaging in Hamiltonian form. Such an approach is used, e.g., in \cite{tresch,  AKN}. In \cite{tresch} bifurcations of the phase portraits of the averaged problem for vertically vibrating suspension point are described. 
 
 We consider the case of arbitrary planar vibrations of the suspension point.  We use the Hamiltonian approach of  \cite{burd} to construct the averaged system and  give a complete description of bifurcations of its phase portraits. Equations for equilibria of the averaged system in this case are obtained in \cite{akul}.

\section{Hamiltonian of the problem}

\begin{figure}[!htbp]
  \centering
  \includegraphics[width=0.35\textwidth]{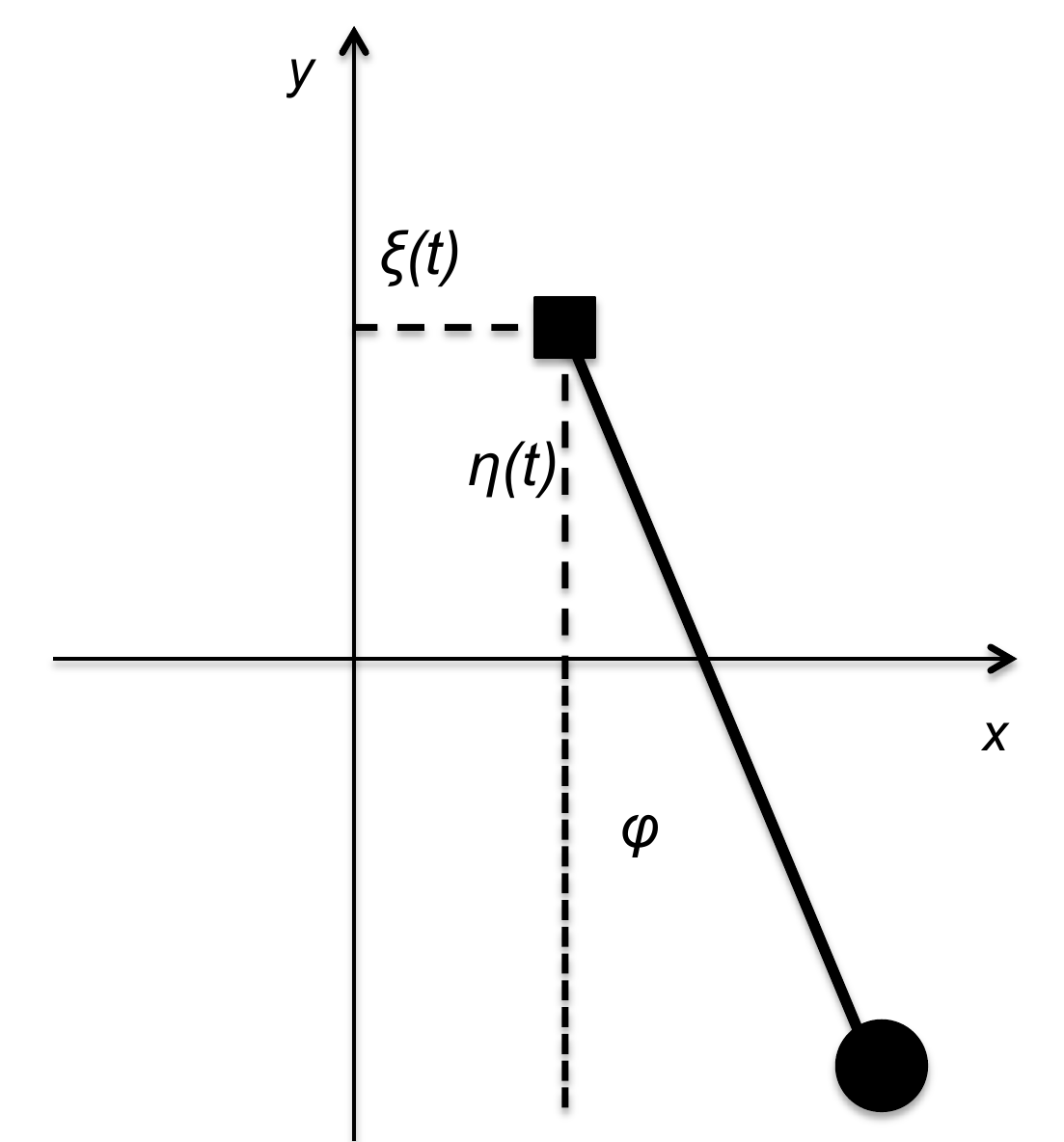}
  \caption{A pendulum with vibrating suspension point.}\label{fig1}
\end{figure}

Consider a simple pendulum, Fig. \ref{fig1}, whose suspension point moves in the vertical plane where the pendulum moves.
Let $l,m$ be length of the massless rod and mass of the bob for this pendulum.  Let $\xi, \eta$ be horizontal and vertical Cartesian coordinates of the suspension point. It is assumed that    $\xi, \eta$ are some given functions of time. Denote by $\ph$ the angle between the pendulum rod and the vertical line straight down. 

Then the kinetic and potential energies of the bob are
\begin{equation*}
T=\frac12\left[l^2\dot \ph^2 +2l \dot \ph(\dot\xi\cos\ph+\dot\eta\sin\ph) +(\dot\xi^2+\dot\eta^2)  \right],\quad V=mg(\eta-l\cos\ph)\, .
\end{equation*} 
Here $g$ is the gravity acceleration.
The generalised momentum conjugate to $\ph$ is
$$
p=\partial T/\partial\dot \ph = ml^2\dot\ph+ml(\dot\xi\cos\ph+\dot\eta\sin\ph)\, .
$$
Thus
$$
\dot\ph = \frac{p}{ml^2}-\frac{(\dot\xi\cos\ph+\dot\eta\sin\ph)}{l}\, .
$$
In order to obtain the Hamiltonian of the problem we should substitute this expression  into $T+V$ and subtract terms which do not depend on $p$. Thus we get the Hamiltonian
$$
H=\frac12\left(\frac{p^2}{ml^2}-2\frac{p(\dot\xi\cos\ph+\dot\eta\sin\ph)}{l} +m(\dot\xi\cos\ph+\dot\eta\sin\ph)^2 \right)-mgl\cos\ph \, .
$$

\section{Averaged Hamiltonian}

Assume that $\xi=\eps\tilde\xi(\omega t/\eps), \eta=\eps\tilde\eta(\omega t/\eps)$, where $\eps$ is a small parameter,  $\tilde\xi(\cdot), \tilde\eta(\cdot)$ are $2\pi$-periodic functions with zero average. Then $\dot \xi, \dot\eta$ are values of order 1. In the system of Hamilton's equations  
$$
\dot \ph = \partial H/\partial p, \quad \dot p = - \partial H/\partial \ph \, .
$$ 
the right hand side is a fast oscillating function of time. In line with the averaging method \cite{bm}, for an approximate description of dynamics of variables $\ph, p$ we average the Hamiltonian $H$ with respect to time $t$. Denote by $A,B,C$ the averages of $\frac{\dot \xi^2}{2}, \frac{\dot \eta^2}{2},$ and $\dot \xi\dot \eta$, respectively. Then the averaged Hamiltonian up to an  additive constant  is
$$
\bar H=\frac12\frac{p^2}{ml^2}+\bar V(\ph), \quad \bar V=\bar V(\ph) =\frac{m}{2}\left[(A-B)\cos 2\ph+C\sin 2\ph\right] - mgl\cos\ph \, .
$$ 
The type of the phase portrait of $H$ is determined by the extrema of $\bar V$. Effectively, it depends on two parameters, $(B-A)/(gl)$ and $C/(gl)$. Our goal is to draw the bifurcation diagram of the problem.  We have to determine a partition of  the plane of parameters by curves into domains such that the type of the phase portrait changes only across these curves. 

The function $\bar V$  is invariant under the transformation $C\to -C, \ph\to -\ph$. Thus, without loss of generality,   in  calculations  below we assume that $C\ge0$.

\section{Particular cases}
\subsection{Case $C=0$}
This is the case considered by Stephenson, Bogolyubov, and Kapitsa under an additional assumption that $A=0$.  Phase portraits for this case are presented in \cite{tresch} (for $A=0$).

\begin{figure}[!htbp]
  \centering
  \includegraphics[width=0.65\textwidth]{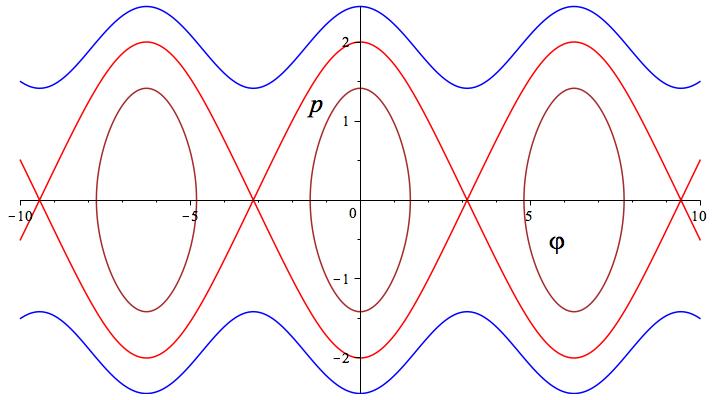}
  \caption{$|\frac{2(B-A)}{gl}|<1 \, .$}\label{fig2}
\end{figure}

If $|\frac{2(B-A)}{gl}|<1$, then $\bar V$ has a minimum at $\ph=0$ and a maximum at  $\ph =\pi$. The phase portrait is shown in Fig. \ref{fig2}  
\\
\\
If $\frac{2(B-A)}{gl}>1 $,  then $\bar V$ has minima at $\ph=0$ and  $\ph=\pi$, and maxima at $\ph = \pm\arccos \frac{gl}{2(A-B)}$.

\begin{figure}[!htbp]
  \centering
  \includegraphics[width=0.65\textwidth]{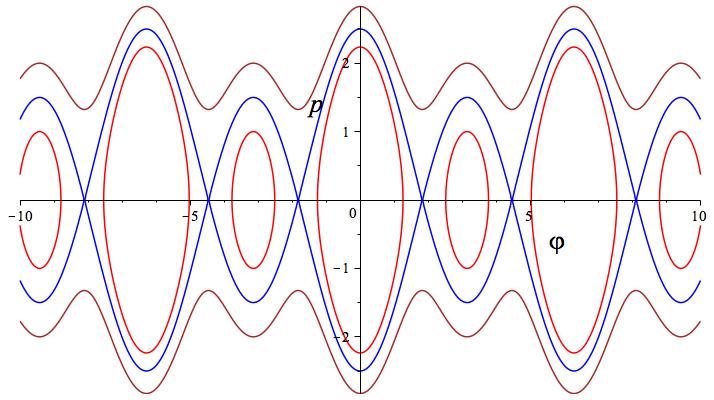}
  \caption{$\frac{2(B-A)}{gl}>1 \, . $}\label{fig3}
\end{figure}

\begin{figure}[!htbp]
  \centering
  \includegraphics[width=0.65\textwidth]{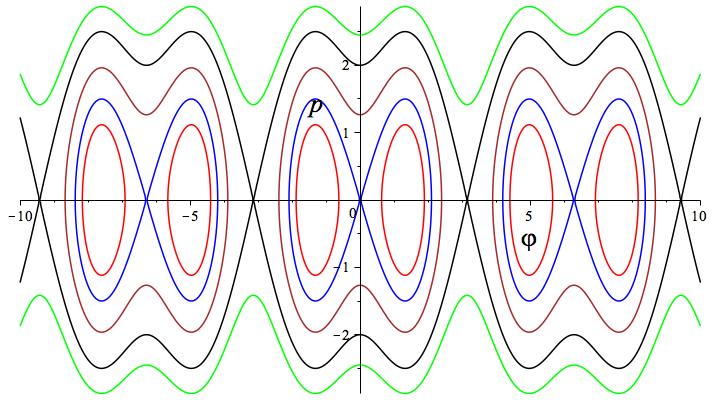}
  \caption{$\frac{2(B-A)}{gl}<-1 \, .$}\label{fig4}
\end{figure}

  The phase portrait is shown in Fig. \ref{fig3}. In this case both the upper and lower  vertical equilibria of the pendulum are stable.
\\
\\
If $\frac{2(B-A)}{gl}<-1 $,  then $\bar V$ has maxima at $\ph =0$ and  $\ph =\pi$, and minima at $\ph =\pm\arccos \frac{gl}{2(A-B)}$.
  The phase portrait is shown in Fig. \ref{fig4}.   In this case both the lower and upper vertical equilibria of the pendulum are unstable.
\\
\\
{ 
We observe a coalescence of two maxima and a minimum of the potential $\bar V$ at  $\frac{2(B-A)}{gl}=1$ and that of two minima  and a maximum of $\bar V$ at $\frac{2(B-A)}{gl}=-1$.

\subsection{Case $A=B$}
The equation for extrema of  $\bar V$ is
$$
\partial \bar V/\partial \ph =  m(C\cos 2\ph + gl\sin\ph)=0 \, 
$$
which implies
\begin{equation}
\label{e_pm}
\sin \ph = \frac{gl\pm\sqrt{g^2l^2+8C^2} } {4C}  \, .
\end{equation}
For the sign ``$-$" in this formula and $C>0$ the value in the right hand side is in the interval $(-1, 0)$. Thus we have two solutions
$$\ph_{-,1}=\arcsin\left( \frac{gl-\sqrt{g^2l^2+8C^2} } {4C} \right), \ \ph_{-,2}= \pi-\ph_{-,1} \, .
$$
They correspond to a minimum and a maximum of $\bar V$, respectively.

If $C> lg$ then for the  sign ``$+$" in (\ref{e_pm}) the value in the right hand side is in the interval \\ (0, 1). In this case  we have two additional solutions
$$
\ph_{+,1}=\arcsin\left( \frac{gl+\sqrt{g^2l^2+8C^2} } {4C} \right), \ \ph_{+,2}= \pi-\ph_{+,1}\, .
$$
They correspond to a maximum and a minimum  of $\bar V$, respectively. 

{ 
Phase portrait for the case $-lg <C<lg$ has the form shown in Fig. \ref{fig2} (but it is in general not symmetric with respect to the axis $\ph=0$, cf. Fig.  \ref{fig9}). 
Phase portraits for the cases  $C>lg$ and  $C<-lg$ are shown in   Fig. \ref{fig6} and Fig. \ref{figx}, respectively.}
\\
\begin{figure}[!htbp]
  \centering
  \includegraphics[width=0.65\textwidth]{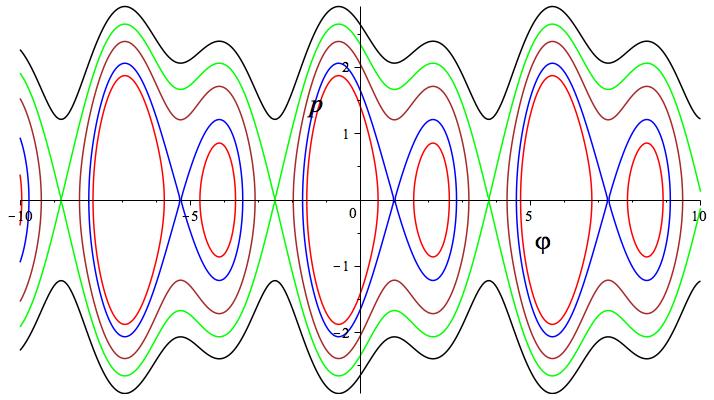}
  \caption{$C>lg \, .$}\label{fig6}
\end{figure}
\\
\\
\begin{figure}[!htbp]
  \centering
  \includegraphics[width=0.65\textwidth]{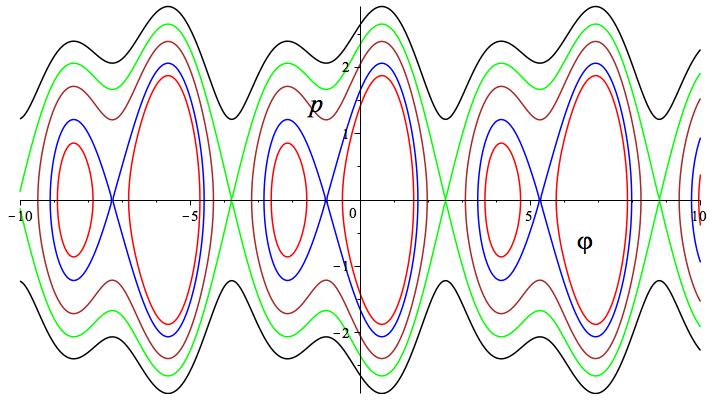}
  \caption{$C<-lg\, .$}\label{figx}
\end{figure}
\\

\section{General  case}
We will build a partition of the parameter plane of the problem into domains
corresponding to different types of phase portraits of the averaged system.
{Boundaries of these  domains are critical curves.} There are two types of
critical curves:

$\bullet$ \  Curves corresponding to degenerate equilibria: first and second
derivatives of $\bar V$ vanish for parameters on these curves. The number of
equilibria changes at crossing  such a curve in the plane of
parameters.

$\bullet$ \  Curves such that the function $\bar V$ has the same value at two different saddle 
equilibria for parameters on these curves. The number of equilibria
remains the same at crossing  such a curve in the plane of
parameters. Curves of the second type {separate   adjacent regions}
 with different behavior of separatrices.

\subsection{Critical curves corresponding to degenerate equilibria}

These curves correspond to values of parameters such that
\begin{eqnarray*}
\partial \bar V/\partial \ph =  m\left[(B-A)\sin 2\ph + C\cos 2\ph + gl\sin\ph \right]=0 \, ,\\
\partial^2 \bar V/\partial \ph^2 =  m\left[2(B-A)\cos 2\ph - 2C\sin 2\ph + gl\cos\ph \right]=0 \, 
\end{eqnarray*}
which is equivalent to
\begin{eqnarray}
\label{deg_0}
\frac{B-A}{gl}\sin 2\ph + \frac{C}{gl}\cos 2\ph + \sin\ph=0 \, , \\
\frac{2(B-A)}{gl}\cos 2\ph - \frac{2C}{gl}\sin 2\ph + \cos\ph=0 \, .\nonumber
\end{eqnarray}
This system is invariant with respect to the change $B-A \to -(B-A), \,  \ph \to \pi-\ph$. Thus, without loss of generality,   in  this subsection we assume that $B-A\ge0$. As agreed earlier, we also assume that $C>0$.

The system (\ref{deg_0}) can be considered as a system of two linear equations for unknown $({B-A})/{(gl)}$ and ${C}/{(gl)}$. Solving it for these unknown we get
\begin{eqnarray*}
\frac{B-A}{gl}=-\sin\ph\sin2\ph-\frac12\cos\ph\cos 2\ph,\\ \frac{C}{gl}=-\sin\ph\cos 2\ph+\frac12\cos\ph\sin 2\ph \, .
\end{eqnarray*} 
An a bit more compact form of these relations is 
\begin{eqnarray}
\label{param}
\frac{B-A}{gl}&=&\cos^3\ph -\frac32\cos\ph ,\\
 \frac{C}{gl}&=&\sin^3\ph \, . \nonumber
\end{eqnarray} 
These relations give a parametric representation with the parameter $\ph$  of the critical curve $\Gamma$ corresponding to degenerate equilibria. This curve is shown in Fig. \ref{fig7}.
\\
\begin{figure}[!htbp]
  \centering
  \includegraphics[width=0.65\textwidth]{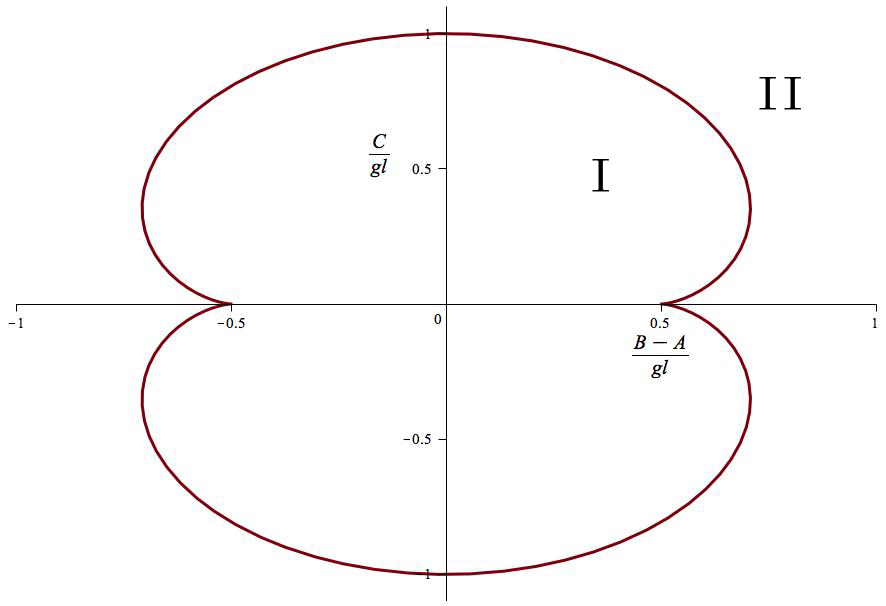}
  \caption{Bifurcation curve corresponding to degenerate equilibria.}\label{fig7}
\end{figure}
\\
The  pendulum has 2 equilibria  if parameters are in domain $I$ (Fig. \ref{fig9}), and 4 equilibria  if parameters are in domain $II$ (Figs. \ref{fig8.1},  \ref{fig8.2}).
\\
\begin{figure}[!htbp]
  \centering
  \includegraphics[width=0.65\textwidth]{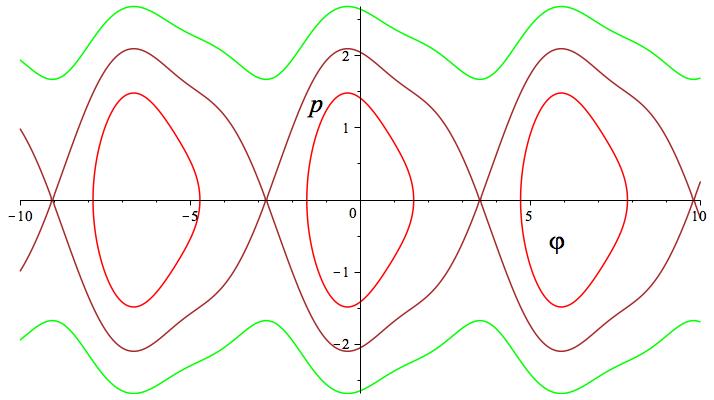}
  \caption{Phase portrait for domain $I$.}\label{fig9}
\end{figure}
\\
\\
\begin{figure}[!htbp]
  \centering
  \includegraphics[width=0.65\textwidth]{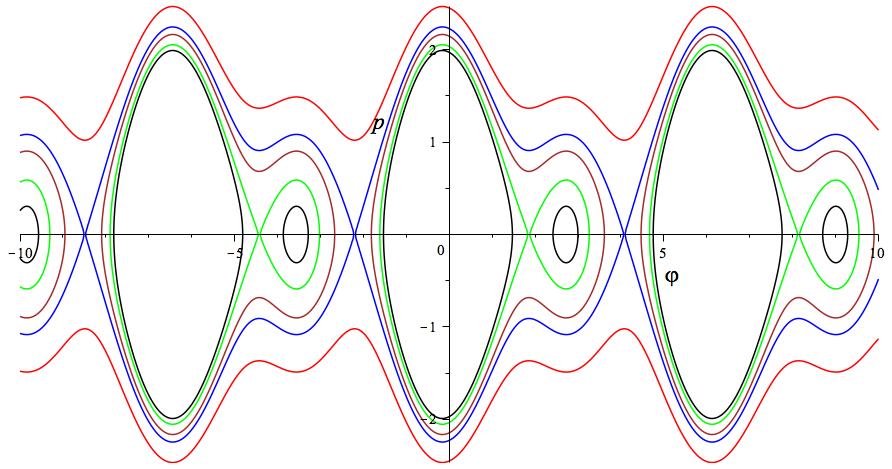}
  \caption{Phase portrait for domain $II,\, C>0.$}\label{fig8.1}
\end{figure}
\\
\\
\begin{figure}[!htbp]
  \centering
  \includegraphics[width=0.65\textwidth]{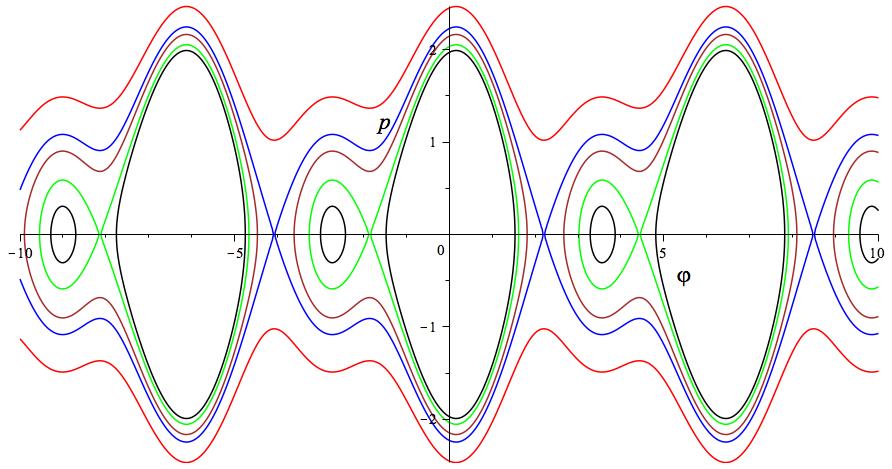}
  \caption{Phase portrait for domain $II, \, C<0.$}\label{fig8.2}
\end{figure}
\\

{ 

The curve $\Gamma$  has two cusp points. Regular points of $\Gamma$ correspond to centre-saddle bifurcations: centre and saddle equilibria coalesce and disappear.    Let us prove that these bifurcations are not degenerate,  i.e., at regular points of $\Gamma$,     (a) $ \partial^3 \bar V/\partial \ph^3 \ne 0$ and (b) the derivative of $ \partial \bar V/\partial \ph  $ along a  transversal to $\Gamma$ direction is different from 0 (after calculation of derivatives one should substitute the equilibrium  value of $\ph$ from  (\ref{deg_0}) into  obtained expressions). 

For (a), we have 
\begin{eqnarray*}
\frac{1}{mgl}\frac {\partial^3 \bar V}{\partial \ph^3} = -\frac{4(B-A)}{gl}\sin 2\ph - \frac{4C}{gl}\cos 2\ph - \sin\ph =3 \sin \ph \ne 0
\end{eqnarray*} 
at regular points of $\Gamma$ (we used  the first equality in (\ref {deg_0}) here).

For (b), we note that   the parametric representation (\ref {param})  of $\Gamma$ implies that the vector
$$
{\bf u} =(-3\cos^2\ph \sin \ph +\frac32\sin\ph ,\ 3\sin^2\ph \cos \ph)\,.
$$
is tangent to $\Gamma$. Then
$$
{\bf u}^{\bot} =(3\sin^2\ph \cos \ph, \ 3\cos^2\ph \sin \ph -\frac32\sin\ph)\,.
$$
is a vector transversal  to $\Gamma$. For the  directional derivative of $ \partial \bar V/\partial \ph  $ along  ${\bf u}^{\bot}$ we get  
\begin{eqnarray*}
\frac{1}{mgl} \nabla_{{\bf u}^{\bot}}  \left(\frac{\partial \bar V}{\partial \ph}\right) &=&
3\sin^2\ph \cos \ph \sin 2\ph + (3\cos^2\ph \sin \ph -\frac32\sin\ph)\cos 2\ph\\
&=&\frac32\sin\ph \left(2\sin\ph \cos \ph \sin 2\ph + (2\cos^2\ph  -1)\cos 2\ph\right)=\frac32\sin\ph\ne 0
\end{eqnarray*}
at regular points of $\Gamma$. Thus, bifurcations at regular points of $\Gamma$ are not degenerate. 

At the cusp points of $\Gamma$ we have  $\partial \bar V/\partial \ph = \partial^2 \bar V/\partial \ph^2=\partial^3 \bar V/\partial \ph^3 =0,\ \partial^4 \bar V/\partial \ph^4= \pm 3mgl \ne 0$ at  the equilibria  ($\pm$ for the left and right cusp points, respectively).  
Near cusp points the curve $\Gamma$ has the standard form of semi-cubic parabola: 
$$
\left\{
\begin{aligned}
\frac{B-A}{gl}=-\frac 12 - \frac34 \ph^2 +\ldots   \\
\frac{C}{gl}=\ph^3+\ldots\phantom{*******}   \\
\end{aligned}
\right.
 \quad {\rm  or} \quad  \quad 
\left\{
\begin{aligned}
\frac{B-A}{gl}\phantom{*}=\frac 12 + \frac34 (\ph-\pi)^2 +\ldots   \\
\frac{C}{gl}=-(\ph-\pi)^3+\ldots\phantom{*******}   \\
\end{aligned}
\right.
$$
for the left and the right cusp points, respectively. This implies that the normal forms of the Hamiltonian near the cusp points are of a generic form presented, e.g., in \cite {HH}, example 2.6, and thus  the bifurcations at  cusp points are not degenerate.}

\subsection{Critical curves corresponding to {saddle} equilibria with equal values of potential energy}

Let $\ph_1$ and $\ph_2$ be coordinates of such equilibria, $\ph_1\ne \ph_2 \ {\rm mod}\, 2\pi$. We should have
$$\partial \bar V(\ph_1)/\partial \ph=0,\, \partial \bar V(\ph_2)/\partial \ph=0, \, \bar V(\ph_1)=\bar V(\ph_2) \, 
$$
or, in the explicit form,
\begin{eqnarray}
\label{three_eq}
\frac{B-A}{gl}\sin 2\ph_1& +& \frac{C}{gl}\cos 2\ph_1 + \sin\ph_1 =0 \, ,\\
\frac{B-A}{gl}\sin 2\ph_2 &+ & \frac{C}{gl}\cos 2\ph_2 + \sin\ph_2 =0 \, , \nonumber \\
-\frac{B-A}{gl}(\cos 2\ph_1-\cos 2\ph_2)&+&\frac{C}{gl}(\sin 2\ph_1 - \sin 2\ph_2) -2(\cos\ph_1-\cos\ph_2) =0 \, . \nonumber
\end{eqnarray}
As  $\ph_1\ne \ph_2 \ {\rm mod}\, 2\pi$, we can divide the last equation by $4\sin\frac{ \ph_1-\ph_2}{2}$. It reduces to
\begin{eqnarray}
\label{four_eq}
\frac{B-A}{gl} \sin(\ph_1+\ph_2) \cos\frac{ \ph_1-\ph_2}{2}+\frac{C}{gl} \cos\frac{ \ph_1-\ph_2}{2}\cos(\ph_1+\ph_2)   +\sin\frac{ \ph_1+\ph_2}{2} =0 \, .
\end{eqnarray}
The first two equations in  (\ref{three_eq}) can be considered as a system of two linear equations for unknown $({B-A})/{(gl)}$ and ${C}/{(gl)}$. 
The determinant of this system is 
$$
D=\sin2(\ph_1-\ph_2) \, .
$$ 
We should consider two cases: $D=0$ and $D\ne 0$.

In the case $D\ne 0$ we solve first two equations in  (\ref{three_eq}) for $({B-A})/{(gl)}$ and ${C}/{(gl)}$ and get
\begin{eqnarray*}
\frac{B-A}{gl}&=&\frac{-\sin\ph_1\cos 2\ph_2+\sin\ph_2\cos 2\ph_1}{\sin2(\ph_1-\ph_2)} \, ,\\
\frac{C}{gl}&=&\frac{\sin\ph_1\sin 2\ph_2-\sin\ph_2\sin 2\ph_1}{\sin2(\ph_1-\ph_2)} \, .
\end{eqnarray*}
These relations can be rewritten in the form
\begin{eqnarray*}
\frac{B-A}{gl}&=&\frac{2\sin\frac{ \ph_1-\ph_2}{2}\cos\frac{ \ph_1+\ph_2}{2}\left( \cos( \ph_1+\ph_2)-\cos( \ph_1-\ph_2)-1       \right)}{\sin2(\ph_1-\ph_2)}  \, ,\\
\frac{C}{gl}&=&\frac{ 2 \sin\frac{ \ph_1-\ph_2}{2}\sin \frac{ \ph_1+\ph_2}{2} \left(-\cos( \ph_1+\ph_2) + \cos( \ph_1-\ph_2)     \right)}    {\sin2(\ph_1-\ph_2)} \, 
\end{eqnarray*}
or, after a cancellation, in the form
\begin{eqnarray*}
\frac{B-A}{gl}&=&\frac{\cos\frac{ \ph_1+\ph_2}{2}\left( \cos( \ph_1+\ph_2)-\cos( \ph_1-\ph_2)-1       \right)}{2\cos( \ph_1-\ph_2)\cos\frac{ \ph_1-\ph_2}{2}} \, ,\\
\frac{C}{gl}&=&\frac{  \sin \frac{ \ph_1+\ph_2}{2} \left(-\cos( \ph_1+\ph_2) + \cos( \ph_1-\ph_2)     \right)} {2\cos( \ph_1-\ph_2)\cos\frac{ \ph_1-\ph_2}{2}} \, .
\end{eqnarray*}
Now we substitute these relations into equation (\ref{four_eq}), and, assuming that $ \ph_1\neq -  \ph_2\ {\rm mod}\,  2\pi$, divide the obtained equation by   $\sin \frac{ \ph_1+\ph_2}{2}$. We get
\begin{eqnarray*}
&\phantom{*}&\frac{\cos^2\frac{ \ph_1+\ph_2}{2}\left( \cos( \ph_1+\ph_2)-\cos( \ph_1-\ph_2)-1       \right)}{\cos( \ph_1-\ph_2)}\\
 &+&\frac{  \left(-\cos( \ph_1+\ph_2) + \cos( \ph_1-\ph_2)     \right)} {2\cos( \ph_1-\ph_2)} \cos(\ph_1+\ph_2)   + 1 =0 \, 
\end{eqnarray*}
which is equivalent to
\begin{eqnarray*}
&\phantom{*}&(1+\cos( \ph_1+\ph_2))\left( \cos( \ph_1+\ph_2)-\cos (\ph_1-\ph_2)-1       \right) \\
 &+&{ \left(-\cos( \ph_1+\ph_2) + \cos( \ph_1-\ph_2)     \right)}  \cos(\ph_1+\ph_2)   +  2\cos( \ph_1-\ph_2)=0 \, 
\end{eqnarray*}
and then equivalent to
$$ \cos( \ph_1-\ph_2)=1 \, 
$$
which can not  be satisfied for  $\ph_1 \neq  \ph_2\, {\rm mod}\,  2\pi$. 

Thus, the only remaining possibilities are related to the special cases  $D=\sin2(\ph_1-\ph_2)=0$ and $ \ph_1= -  \ph_2\, {\rm mod}\,  2\pi$. It is easy to check that each of these relations together with relations (\ref{three_eq}) implies $C=0$. As $\ph_1$ and $\ph_2$ should correspond to different  saddles, we should consider the domain were there are four equilibria. Thus, {possible} bifurcational curves are rays $\{C=0, \frac{2(B-A)}{gl}>1\} $  and $\{C=0, \frac{2(B-A)}{gl}<-1 \}$. {The equilibria with equal values of potential energy are saddles for the first  and centres for the second of these rays (see Figs.  \ref{fig3} and  \ref{fig4}).  Thus the only bifurcational curve of the considered type is the ray $\{C=0, \frac{2(B-A)}{gl}>1\} $.}


\subsection{Bifurcation diagram}

The partition of the plane of parameters into domains with qualitatively different phase portraits is shown in Fig. \ref {fig8}. The  pendulum has 2 equilibria  if parameters are in domain $I$ (phase portrait in Fig. \ref{fig9}), and 4 equilibria  if parameters are {in domain  $II$ (phase portrait in Fig. \ref{fig8.1} for $C>0$ 
 and Fig.  \ref{fig8.2} for $C<0$; these phase portraits are topologically the same but have different types of asymmetry of the eight-shaped separatrix).   The boundary $\Gamma$ between these domains is  the curve  of degenerate equilibria.  The ray $\{C=0, \frac{2(B-A)}{gl}>1\} $  corresponds to phase portraits with a heteroclinic trajectory (Fig. \ref{fig3}).   Passage through this  ray transforms  the phase portrait of Fig. \ref{fig8.1} via that of Fig. \ref{fig3} to that of   Fig. \ref{fig8.2}.     At approach the curve  of degenerate equilibria from  the domain  $II$ at $C>0$ (respectively $C<0$),} the right (respectively, the left)  loop of eight-shaped separatrix shrinks to a point and disappears.  
\\
\\
\begin{figure}[!htbp]
  \centering
  \includegraphics[width=0.65\textwidth]{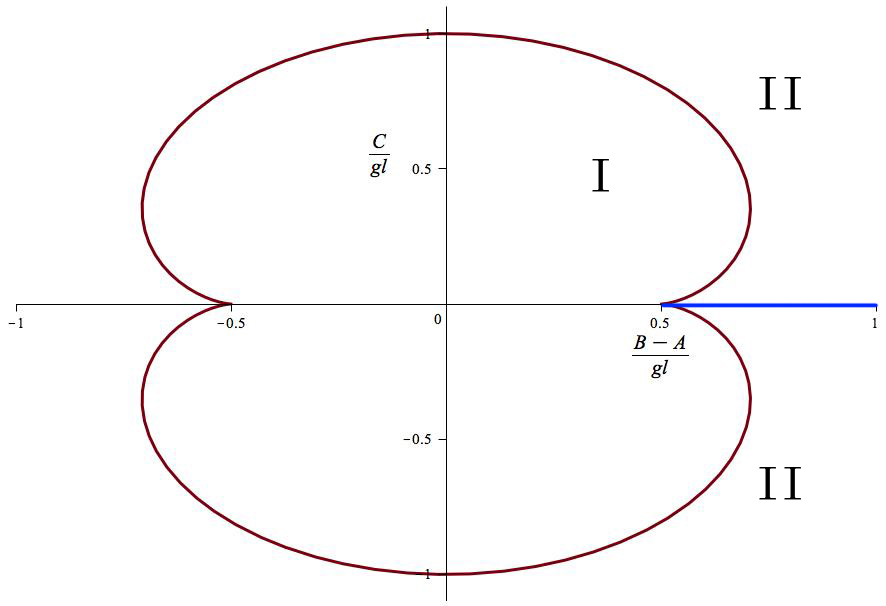}
  \caption{Partition of the parameter plane.}\label{fig8}
\end{figure}
\\

\section{Conclusion}
Our study in this paper  provides a complete description of bifurcations of phase portraits for a pendulum with vibrating suspension point in the approximation of the averaging method. The relation of such a description to dynamics in the exact (not averaged) problem for high-frequency oscillations is discussed, e.g., in \cite{AKN}, Sect. 6.3.3.B.  One should consider  the Poincar\'e return map for the plane $\ph, p$ (Poincar\'e section $t=0 \ {\rm mod}\, \frac{2\pi\eps}{\omega}$). Fixed points of this map are close to equilibria of the averaged system. The plane $\ph, p$ is filled by invariant curves of this map up to a reminder of a measure which is exponentially small in the perturbation parameter  $\eps$.  Stable and unstable manifolds of hyperbolic fixed points of the return map are split,  but this splitting is exponentially small. These manifolds are close to separatrices of the averaged system. Thus, the study of the averaged system  provides  considerable information about  dynamics of the exact (not averaged) problem. 

\medskip
{ 
{\bf Acknowledgement}

The authors are very grateful to the anonymous referee for a constructive criticism and useful suggestions.}

\end{document}